\documentclass[11pt]{article}
\usepackage{amssymb}

  \newtheorem {Theorem} {Theorem} [section]
  
  \newtheorem {Proposition} [Theorem] {Proposition}
  \newtheorem {Lemma} [Theorem] {Lemma}

  \newtheorem {rem} [Theorem] {}

  \newenvironment{Proof}[1][Proof] 
  {       
    \emph{#1:} 
    \setlength{\parskip}{0ex plus0ex minus0ex}}
  { \hspace*{\fill} \ensuremath{\square}   }

    \newcommand {\betrag} [1] {\ensuremath{ \left\vert  #1  \right\vert } } 
    \newcommand {\norm} [2] [] {\ensuremath{ \left\Vert  #2  \right\Vert_{#1} } } 
    \newcommand {\R} {\ensuremath{\mathbb{R}}}
    
    \newcommand {\N} {\ensuremath{\mathbb{N}}}
    
    \newcommand {\F}  {\ensuremath {\mathcal{F}}}
    \newcommand {\Om} {\ensuremath{\Omega}}

    \newcommand {\T} {\mathcal{T}}

    \newcommand {\spec} {\mathrm{spec}}
    
    \newcommand {\W} {\mathcal{W}}

    \newcommand {\x} {\mathbf{x}}
    \newcommand {\y} {\mathbf{y}}
    \newcommand {\z} {\mathbf{z}}
    \newcommand {\muT} {\mu_{\scriptscriptstyle T}}
    
    \newcommand {\cH} {\mathcal{H}}

\begin{document} 

\title
{\Large \bf Existence of Gibbs measures relative to Brownian motion}
\author{
\small Volker Betz \\[0.1cm]
{\it \small Zentrum Mathematik, Technische Universit\"at M\"unchen} \\
{\it \small Gabelsbergerstr. 49, 80290 M\"unchen, Germany} \\
{\small betz@mathematik.tu-muenchen.de}\\ [1cm]}
\date{}
\maketitle
\begin{abstract}
We prove existence of infinite volume 
Gibbs measures relative to Brownian motion. We require the pair potential $W$ to
fulfill a uniform integrability condition, but otherwise our restrictions on
the potentials are relatively weak.  In particular,
our results are applicable to the massless Nelson model. 
We also prove an upper bound for path fluctuations under the infinite volume Gibbs measures.
\end{abstract}

\noindent {\small
{\em Keywords: } Stationary non-Markov processes, Gibbs measures, Nelson model.\\
{\em Mathematics subject classification: } 60G10, 82B99 \\}

\vspace{0.5cm}

\section{Introduction} \label{S1}

Let us define a probability measure on  $C([-T,T],\R^d)$ by 
\begin{equation} \label{fvgibbs with bc}
d\mu_T^{y,z} (x)  = \frac{1}{Z_T(y,z)} e^{ - \int_{-T}^T V(x_s) \, ds - 
 \int_{-T}^T ds \int_{-T}^T dt \, W(x_t,x_s,|t-s|)} \, d\W^{y,z}_{[-T,T]}(x).
\end{equation}
on $C([-T,T],\R^d)$.
Here,  $T>0$, $y,z \in \R^d$,  
$\W^{y,z}_{[-T,T]}$ is pinned Brownian motion starting in $y$ at time $-T$ and
ending in $z$ at time $T$, $Z_T(y,z)$ normalizes $\mu_T^{y,z}$ to a probability 
measure, and
$V: \R^d \to \R$ and $W: \R^d \times \R^d \times \R \to \R$ are measurable functions
with some additional properties to be specified later. One choice of $V$ and $W$ that
fits into the framework of the present paper is 
\begin{equation} \label{VWexample}
d=3,\qquad  V(x)= - 1/ |x| \quad \mbox{ and }\quad  W(x,y,t) = - 1/(|x-y|^2 + t^2 + 1).
\end{equation}

A natural problem in the context of (\ref{fvgibbs with bc}) is the existence and uniqueness,
i.e. independence of the `boundary conditions' $y,z \in \R^d$, 
of a limiting probability
measure $\mu_{\infty}$ on $C(\R,\R^d)$ as $T \to \infty$ in (\ref{fvgibbs with bc}).
$\mu_{\infty}$ will be called (infinite volume) Gibbs measure relative to Brownian motion; 
this terminology already suggests  a close relationship with statistical mechanics.
We will outline this connection as well as a link to the theory of large deviations 
toward the end end of this introduction, but first let us study (\ref{fvgibbs with bc}) 
in its own right. 

An easy special case of (\ref{fvgibbs with bc}) is obtained by choosing $W=0$. Then, via the 
Feynman-Kac formula, $\mu_T^{x,y}$ is related to the Schr\"odinger operator $H_0 = -\frac{1}{2}
\Delta + V$. If 
\begin{equation} \label{mincond V}
\mbox{$H_0$ has a ground state $\psi_0 \in L^2(\R^d)$,}
\end{equation}
then the infinite volume Gibbs measure exists and is given by the stationary solution of the 
stochastic differential equation $dX_t = \frac{\nabla \psi_0}{\psi_0}(X_t) + dB_t$ 
(see  \cite{Si79} or equations (\ref{SDE}) and (\ref{FKF})).
We will take the point of view that the $W \neq 0$ case is a perturbation of the $W=0$ case.
The existence problem for $\mu_{\infty}$ 
can then be regarded as a generalization to the problem of  finding
stationary solutions for stochastic differential equations. An important difference of the two 
problems is that, unlike solutions to stochastic differential equations, the limiting measure
will not be the measure of a Markov process if $W \neq 0$.

When looking for reasonable conditions on $V$ and $W$ that ensure existence of $\mu_{\infty}$, 
a natural requirement on $V$ is that it should lead to an infinite volume Gibbs measure at
least in the case $W=0$. (\ref{mincond V}) is a sufficient condition for this. 
As far as the `perturbation' $W$ is concerned, we should require that its effect does not completely
outweigh the effect of the $V$. In other words, $W$ has to be extensive, i.e.
\begin{equation} \label{extensive}
\limsup_{T \to \infty} \frac{1}{T} \left| \int_{-T}^T ds \int_{-T}^T dt \, W(x_s,x_t,|s-t|) \right| < \infty
\end{equation} 
at least for a reasonable class of $x \in C(\R,\R^d)$. While 
(\ref{extensive}) may not be sufficient for the existence of $\mu_{\infty}$ in general, 
additional conditions on $W$ should be more of a technical nature.

As for uniqueness, already the case $W=0$ shows \cite{BL00} that we can only 
expect $\mu_{\infty}$ to be unique among the measures supported on
a subset of $C(\R,\R^d)$ which is characterized by a condition on the growth of paths at infinity. 
Once this restriction is made, according to the folklore a sufficient condition is 
that the interaction energy 
\begin{equation} \label{interaction}
 I = \sup_{x \in C(\R, \R^d)} \left| \int_{-\infty}^0 ds \int_0^{\infty} dt \,  W(x_t,x_s,|t-s|) \right| 
\end{equation}
between left and right half of the path is finite. Such a strong result is not available at present, but 
\cite{OS00} and \cite{LM00} have some results about uniqueness, and \cite{OS00} gives an
example where uniqueness fails when (\ref{interaction}) is not fulfilled. In the present work, we 
have nothing to say about uniqueness, focussing on existence instead.

Several authors have by now studied the existence problem. All of them 
assume (\ref{extensive}) in some form, but also need additional restrictions on $V$ and $W$.
In \cite{OS00}, the first mathematical account on the subject, 
correlation inequalities are used, and consequently the potentials $V$ and $W$ have to  
fulfill certain convexity assumptions. In \cite{LM00}, a cluster expansion method is applied,
requiring a small parameter (coupling constant) in front of 
$W$ as well as a $V$ that is growing faster than quadratically at infinity. Recently, \cite{Ha01}
used an integration by parts formula. His restrictions on $W$ are weak, but strong assumptions 
on the asymptotic behaviour of $V$ are needed. In particular, $V$ has to grow at least quadratically 
at infinity.

In this work we establish a new method for proving existence of $\mu_{\infty}$, relying on a stopping 
time estimate. The main advantage over the existing approaches is that our restrictions on $V$ are
almost as weak as (\ref{mincond V}). All cases from \cite{OS00,LM00,Ha01} are covered, and in addition 
we allow for $V$'s which do not grow at infinity. For the pair potential $W$, the main assumption 
essentially is that (\ref{extensive}) holds uniformly in $x \in C(\R,\R^d)$. In addition, we need 
a  `pathwise shift condition' that is somewhat implicit but easy to verify for many concrete examples of $W$.
If we assume that $V$ fits in the framework of \cite{OS00}, \cite{LM00} or \cite{Ha01}, then 
on the one hand the (uniform) integrability conditions on $W$ that we 
impose are stronger than those needed there. On the other hand,
we neither need the convexity assumed in \cite{OS00}, nor the small parameter of \cite{LM00}, nor
the differentiability needed in \cite{Ha01}. An important feature that our work shares with all of 
the above is that the interaction energy (\ref{interaction})
between the left and the right half-line is not assumed to be finite. 

As mentioned before, there exist connections or (\ref{fvgibbs with bc}) with statistical mechanics
as well as with the theory of large deviations. The latter connection is seen most clearly when we
replace the exponent in (\ref{fvgibbs with bc}) by
\begin{equation} \label{mean field}
-\int_{-T}^T V(x_s) \, ds - \frac{1}{2T} \int_{-T}^T ds \int_{-T}^T dt \, \tilde{W}(x_s,x_t) 
\end{equation}
with some nice function $\tilde{W}$. (\ref{mean field}) is then a functional of the local time,
and thus is a special case of the theory of Donsker and Varadhan \cite{DoVa76}. 
So in a sense, these systems are extremely well understood. It turns out that the limiting process
for interactions like (\ref{mean field}) is a Markov process. This is not the case 
for the actual system (\ref{fvgibbs with bc}), which shows that although (\ref{mean field}) and
(\ref{fvgibbs with bc}) may look similar, they yield very different limiting objects.
In the language of statistical mechanics,
(\ref{mean field}) is a mean field interaction, while  (\ref{fvgibbs with bc}) is a local interaction.

To link (\ref{fvgibbs with bc}) with statistical mechanics, more precisely with the theory of 
lattice spin systems, we 
discretize (\ref{fvgibbs with bc}) by replacing Brownian motion with a random
walk with state space $\R^d$ and Gaussian step size distribution. We then obtain
a finite volume Gibbs measure on a one-dimensional system of $\R^d$-valued spins, with 
single site potential $V$, quadratic nearest neighbour interaction and long range pair interaction $W$.
The reference measure is the product of $d$-dimensional Lebesgue-measures. An equivalent description
of this spin system, a little bit closer to  (\ref{fvgibbs with bc}), is to incorporate the
nearest-neighbour interaction into the reference measure, which then becomes the measure of a random 
walk pinned at $-T$ and $T$. 

Although we will not do it here, our method can be easily adapted to the 
lattice context, where it yields a new way of proving existence of Gibbs measures for
one-dimensional systems of unbounded spins. For such systems, extremely powerful methods are already 
available: there is the superstability estimate by D. Ruelle \cite{Ru70}, applied in \cite{LP76}, 
which has the big advantage of not being restricted to one-dimensional systems; there are the results
of R. L. Dobrushin \cite{Do73,Do73a}, which are valid only for one-dimensional systems, but extremely
general otherwise. However, superstability corresponds to rapidly growing single site potential,
while one of Dobrushin's few restrictions is that the interaction energy between left and right
half-space must be bounded. Thus our method covers some new situations in the discrete context also.

Finally, although it should have become clear that Gibbs measures are interesting objects 
also from a purely probabilistic point of view, 
the original motivation for studying them is a physical one. 
Nelson \cite{Ne64} first used measures with a structure similar to 
(\ref{fvgibbs with bc}) with a $W$ of the type given in (\ref{VWexample}) 
to study the ultraviolet divergence in a model of a 
quantum particle coupled to a scalar bosonic field, nowadays known as Nelson's model. 
In \cite{Sp86}, Gibbs measures are used to estimate the effective mass of the polaron.  
Recently \cite{LMS00,BHLMS} study various aspects of the ground state of
Nelson's model by using Gibbs measures.

\section{Finite volume Gibbs measures} \label{S2}

We start by specifying conditions on the potentials $V$ and $W$ appearing in 
(\ref{fvgibbs with bc}). 
A measurable function $V: \R^d \to \R$ is said to be in the Kato class  \cite{Si82}, 
$V \in \mathcal{K}(\R^d)$, if
$$ \sup_{x \in \R} \int_{\{|x-y| \leq 1\}} |V(y)| \, dy < \infty \qquad \mbox{in case $d=1$,} $$
 and
$$ \lim_{r \to 0} \sup_{x \in \R^d} \int_{\{|x-y| \leq r\}} g(x-y) |V(y)| \, dy = 0 \qquad 
\mbox{in case $d \geq 2$.}$$ 
Here,
$$ g(x) = \left\{ \begin{array}{ll} 
        -\ln|x| & \mbox{if } d=2 \\
    |x|^{2-d} & \mbox{if } d \geq 3.
    \end{array} \right. $$
$V$ is locally in the Kato class, $V \in \mathcal{K}_{\mathrm{loc}}(\R^d)$, if
$V 1_K \in \mathcal{K}(\R^d)$ for each compact set $K \subset \R^d$. 
$V$ is Kato-decomposable \cite{BHL00} if 
$$ V = V^+ - V^- \quad \mbox{with} \quad V^- \in \mathcal{K}(\R^d), 
V^+ \in \mathcal{K}_{\mathrm{loc}}(\R^d),$$
where $V^+$ is the positive part and $V^-$ is the negative part of $V$.

Our conditions on $V$ are:
\begin{itemize} 
\item[(V1):] $V: \R^d \to \R$ is Kato-decomposable.
\item[(V2):] The Schr\"odinger operator 
$$ H_0 = -\frac{1}{2} \Delta + V$$ 
(where $\Delta$ denotes the Laplace operator) acting in $L^2(\R^d)$ 
fulfills $\inf\spec(H_0)=0$. Moreover, $H_0$ has a unique, strictly positive ground 
state $\psi_0 \in L^2(\R^d) \cap L^1(\R^d)$, i.e. $0$ is an eigenvalue of multiplicity one
with corresponding eigenfunction $\psi_0$.
\end{itemize}

Condition (V1) guarantees that the factor $\exp(-\int_{-T}^T V(x_s) \, ds)$ appearing
in (\ref{fvgibbs with bc}) is integrable with respect to Brownian motion \cite{Si82}. 
The existence
of a ground state in (V2) ensures the existence of an infinite volume Gibbs measure in case
$W=0$, while $\inf\spec(H_0)=0$ is included for convenience and 
can be achieved by simply adding a constant to 
$V$ and changing the normalizing constant in  (\ref{fvgibbs with bc}) accordingly. 
Finally, $\psi_0 \in L^1$ will be needed in the proof of Theorem \ref{Th3.2}, 
but is only a mild 
restriction, since in most cases of interest $\psi_0(x)$ decays exponentially
for large $x$ \cite{Ca78}. 

Examples for potentials $V$ that fulfill $(V1)$ 
and $(V2)$ are continuous functions bounded below and growing at infinity, as
well as functions bounded above but with the negative part having Coulomb type singularities.

Schr\"odinger operators with Kato-decomposable potentials have many nice properties \cite{Si82}. 
In this paper we will need the fact that the kernel $K_t(x,y)$ of $e^{-tH_0}$ uniformly bounded 
and bounded away from zero on compact sets, and that $y \mapsto K_t(x,y)$ is integrable uniformly 
in $x$.

Turning to conditions on $W$, let us write
\begin{equation} \label{H}
\cH_{\Lambda}(x) = - \int \!\!\! \int_{\Lambda} W(x_t,x_s,|t-s|) \, ds \, dt \qquad \qquad
(x \in C(\R, \R^d))
\end{equation}
with $\Lambda \subset \R^2$. In case $\Lambda = [-T,T]^2$, we simply write $\cH_T(x)$.
$C^{(0)}(\R,\R^d)$ will denote functions
which are continuous with the possible exception of the point $0$ but have
left and right hand side limits there. For $\tau > 0$ consider the map  
\begin{equation} \label{theta}
\theta^{(0)}_{\tau}: C(\R, \R^d) \to C^{(0)}(\R,\R^d), \quad 
(\theta^{(0)}_{\tau} x)_t = \left\{ \begin{array}{ll} x_{t+\tau} & \mbox{if } t \geq 0,\\
                                                x_{t-\tau} & \mbox{if } t < 0.
                              \end{array} \right.
\end{equation}
Finally, put
\begin{equation} \label{alpha}
\alpha = \liminf_{|x| \to \infty} V(x) \leq \infty.
\end{equation}
From the way this constant will enter into our proofs it will be clear
that really the quantity 
$$\liminf_{|x| \to \infty} V(x) - \inf\spec(H_0)$$
is the important one, a fact that is obscured by our choice $\inf\spec(H_0)=0$ in (V2).

Our conditions on $W$ are
\begin{itemize} 
\item[(W1):] There exists $C_{\infty} < \infty$ such that
\begin{equation} \label{Cinfty}
\int_{-\infty}^{\infty} |W(x_0,x_s,|s|)| \, ds < C_{\infty} \quad \mbox{and} \quad 
\int_{-\infty}^{\infty} |W(x_s,x_0,|s|)| \, ds < C_{\infty},
\end{equation}
uniformly in $x \in C(\R,\R^d)$.
\item[(W2):] There exist $D \geq 0$ and $0 \leq C < \alpha$ such that
\begin{equation} \label{cond}
\cH_T(x) \leq \cH_T(\theta^{(0)}_{\tau} x) + C \tau + D
\end{equation} 
for all $T, \tau >0$ and all $x \in C(\R,\R^d)$.
\end{itemize}
An immediate consequence of (W1) is
\begin{equation} \label{estimate}
\betrag{ \cH_{\R \times [-S,S]}(x)} \leq 2 C_{\infty}S, \quad \mbox{and} \quad
\betrag{ \cH_{[-S,S] \times \R}(x)} \leq 2 C_{\infty}S.
\end{equation}
(\ref{estimate}) will be used frequently below. 

(W2) looks a little mysterious at first, but the proof of Theorem 
\ref{Th3.2} will show how it comes about naturally. 
To see when (W2) is fulfilled, note that by (\ref{estimate}),
\begin{eqnarray*}
 - \int_0^T \!\! ds \int_0^T \!\! dt \, W(x_t,x_s,|t-s|) & \leq & 4 C_{\infty} \tau -
\int_{\tau}^{T+\tau} \!\!\! ds \int_{\tau}^{T + \tau} \!\!\! dt \, W(x_t,x_s,|t-s|) = \\
 & = & 4 C_{\infty} \tau - \int_0^T \!\! ds \int_0^T \!\! dt W(x_{t+\tau},x_{s+\tau},|t-s|),
\end{eqnarray*}
and similarly for the region $[-T,0]^2$. Thus, if we suppose 
$$ I = \sup_{x \in C(\R,\R^d)} 
\int_{-\infty}^0 \!\! ds \int_0^{\infty} \!\! dt \,|W(x_t,x_s,|t-s|)| < \infty, $$
then $8 C_{\infty} < \alpha$ is a sufficient condition for (W2).
In case $I = \infty$, it is not hard to see that if there exist $L,M > 0$ with 
\begin{equation} \label{Suffcond}
  \int_{-T}^0 \!\! ds \int_0^{T} \!\! dt \big( W(x_s,x_t,|s-t|) - 
  W(x_s,x_t,|s-t| + 2 \tau) \big) \leq L \tau + M
\end{equation}
uniformly in $x \in C(\R,\R^d)$ and $T>0$,
then $12 C_\infty + L < \alpha$ is a sufficient condition for (W2). 
(\ref{Suffcond}) can be checked directly for many choices of $W$, and is
in particular true if $t \mapsto W(x,y,t)$ is increasing for $t>0$ and each fixed
$x,y \in \R^d$. This covers the physically important case 
$$W(x,y,|t|) = -\frac{1}{(|x-y|^2+|t|^2+1)}$$
of the massless Nelson model \cite{BFS99,LM00}.
On the other hand, for
$$ W(x,y,|t|) = \left\{ \begin{array}{ll} 
             - \frac{1}{|t|^2+1} & \mbox{if } |x-y| \leq 2t \\
             0 & \mbox{otherwise} 
   \end{array} \right.$$
$(x,y \in \R)$ together with the path $x_t = t$, we find that 
$ \int_{-T}^0 ds \int_0^T dt \, W(x_s,x_t,|t-s|)$ diverges as $T \to \infty$, but
e.g. $\int_{-T}^0 ds \int_0^T dt \, W(x_{s-1},x_{t+1},|t-s|) = 0$.
Thus (W2) need not hold in general.

We now construct finite volume Gibbs measures.  We will take a point of view that
differs slightly from the one taken in equation \ref{fvgibbs with bc} by incorporating 
the single site potential $V$ 
into the reference measure. This leads to a $P(\phi)_1$-process \cite{Si79}.
To make the paper reasonably self-contained, we include
a short description of this process.

The $P(\phi)_1$-process corresponding to the potential $V$ is the stationary
solution of the stochastic differential equation 
\begin{equation} \label{SDE} 
 dX_t = \frac{\nabla \psi_0}{\psi_0} (X_t) \, dt + dB_t, 
\end{equation}
where $B_t$ denotes Brownian motion in $\R^d$. Remember that $\psi_0$ is the 
ground state of $H_0$. The measure on $C(\R,\R^d)$
corresponding to this process will be denoted by $\mu_0$ and identified with the 
process. $\mu_0$ is a stationary strong Markov process with generator 
$\tilde{H}_0 = \psi_0^{-1} H_0 \psi_0$, where $\psi_0$
and $\psi_0^{-1}$ denote operators of multiplication.

The tool that links (\ref{SDE}) and (\ref{fvgibbs with bc}) is the Feynman-Kac formula.
It says that for a bounded interval $I = [0,T] \subset \R$ and a
$\mu_0$-integrable, $\F_I$-measurable function $f$,
\begin{equation} \label{FKF}
\int f(x) \, d\mu_0(x) = \int \psi_0(x_0) e^{- \int_0^T V(x_s) \, ds} f(x) \psi_0(x_T) 
\, d\W(x).
\end{equation}
Here, $\W$ denotes the infinite mass Wiener measure, and $\F_I$ is the $\sigma$-field 
over $C(\R,\R^d)$ generated by the point evaluations with points inside $I$. A corresponding 
notation for $\sigma$-fields will be used throughout the paper.

By (\ref{FKF}),  the 
invariant measure of $\mu_0$ has the Lebesgue-density $\psi_0^2$. Moreover,
a refined version of the Feynman-Kac formula \cite{Si79} shows that 
the transition density of $\mu_0$ is given in terms of the kernel $K_t(x,y)$ of $e^{-tH}$ by
\begin{equation} \label{Transdens}
E_{\mu_0} (f(x_t) | \F_{\{0\}})(y) = \frac{1}{\psi_0(y)} \int K_t(y,z) \psi_0(z) f(z) \, dz 
\qquad (y \in \R^d).
\end{equation}

We perturb the process $\mu_0$ by the pair potential $W$, i.e. for $T>0$ we 
define the probability measure $\muT$ on $C(\R,\R^d)$ by  
\begin{equation} \label{muT}
d\muT(x) = \frac{1}{Z_T} e^{- \cH_T(x)} \, d\mu_0(x),
\end{equation}
where 
$$ Z_T = \int e^{ - \cH_T(x)} \, d\mu_0(x)$$
is the normalizing constant. Comparing (\ref{muT}) and (\ref{fvgibbs with bc}), we 
see that
instead of pinning the path to $y,z \in \R$ at time $-T$ resp. $T$ (``sharp
boundary condition''), we now allow it to fluctuate according to $\mu_0$ outside 
$[-T,T]$, resulting in a ``smeared-out boundary condition''. This is technically easier
to handle and, as we will see in the course of the paper, good enough to prove  
existence of an infinite volume Gibbs measure.

Let us now check that 
$\muT$ is a finite volume Gibbs measure with respect to the potential 
$W$ and the reference measure $\mu_0$.  
For $S > 0$ write $\T_S$ instead of $\F_{[-S,S]^c}$, and for $\bar{x} \in C(\R,\R^d)$
denote by $\mu_{0}^{S,\bar{x}}$
the version of the regular conditional expectation
$\mu_0(.|\T_S)$ that is given by 
\begin{equation} \label{mu0x}
d\mu_{0}^{S,\bar{x}}(x) = \frac{1}{Z^S(\bar{x})} 
\exp \left(- \int_{-S}^S V(x_s) \, ds \right) \, 
d \left(\W_{[-S,S]}^{\bar{x}} \otimes \delta_{[-S,S]^c}^{\bar{x}}\right)(x).
\end{equation}
Here, $\delta_{[-S,S]^c}^{\bar{x}}$ is the point measure on $C([-S,S]^c ,\R^d)$ 
concentrated in $\bar{x}|_{[-S,S]^c}$,  $\W_S^{\bar{x}}$ is pinned Brownian motion
starting at time $-S$ in $\bar{x}(-S)$ 
and ending at time $S$ in $\bar{x}(S)$, and $Z^S(\bar{x})$ is the normalizing constant.
Moreover, for $S < T$ define  
\begin{eqnarray}
\Lambda(S,T) & = & ([-T,T] \times [-S,S]) \cup ([-S,S] \times [-T,T]) \subset \R^2, \quad \mbox{and}
  \label{LambdaST} \\
d\muT^{S,\bar{x}}(x) & = & \frac{1}{Z_T^{S}(\bar{x})} 
\exp ( \cH_{\Lambda(S,T)}(x)) \, d\mu_0^{S,\bar{x}}(x). \label{finvolgibbs}
\end{eqnarray}
In (\ref{finvolgibbs}), $Z_T^S(\bar{x}) = E_{\mu_0^{S,\bar{x}}}(e^{\cH_{\Lambda(S,T)}})$ 
is again the normalizing constant.

\begin{Lemma}  \label{fvdlr}
For each $S<T$, $\bar{x} \mapsto \muT^{S,\bar{x}}$ is a version of the regular 
conditional expectation $\muT(. | \T_S)$. In other words, $\muT$ is a 
(finite volume) Gibbs measure 
with reference measure $\mu_0$ and potential $W$.
\end{Lemma}

\begin{Proof}
Let $f,g \in L^{\infty}(C(\R,\R^d))$, and suppose $g$ is $\T$-measurable. Then 
\begin{eqnarray*}
\lefteqn{ Z_T \int g(\bar{x}) E_{\muT^{S,\bar{x}}}(f) \, d\muT(\bar{x}) = }\\
 & = & E_{\mu_0} \left( 
       E_{\mu_0} \Big( g \frac{1}{E_{\mu_0}(e^{ \cH_{\Lambda(S,T)}} | \T_S)} 
       E_{\mu_0}(f e^{ \cH_{\Lambda(S,T)}}|\T_S) e^{\cH_{[-T,T]^2}} \Big| \T_S \Big)
       \right) = \\
 & = & E_{\mu_0} \left( 
        g \frac{ E_{\mu_0}(f e^{ \cH_{\Lambda(S,T)}}|\T_S)}
        {E_{\mu_0}(e^{ \cH_{\Lambda(S,T)}} | \T_S)} 
       e^{\cH_{[-T,T]^2 \setminus \Lambda(S,T)}}  
       E_{\mu_0} \Big( e^{\cH_{\Lambda(S,T)}} \Big| \T_S\Big)
       \right) = \\
 & = & Z_T E_{\muT}(fg).
\end{eqnarray*}
Dividing by $Z_T$ finishes the proof.
\end{Proof}

\section{Infinite volume Gibbs measures} \label{S3}

We say that a sequence $\nu_n$ 
of measures on $C(\R,\R^d)$  converges locally weakly to a measure $\nu$ 
if for each bounded interval $I \subset \R$, the restrictions of $\nu_n$ to 
$\F_I$ converge weakly to the restriction of $\nu$ to $\F_I$.
It is easy to see that, when $C(\R,\R^d)$ equipped with the topology of uniform 
convergence on compact sets, local weak convergence is equivalent to weak convergence.

An infinite volume cluster point of a family $(\nu_T)_{T>0}$ of probability measures  
is a cluster point of any sequence $(\nu_{t_n})_{n \in \N}$, where $t_n \to \infty$ 
as $n \to \infty$. 

We will show that the family $(\muT)$ is relatively compact in the topology of 
local weak convergence. From this the existence of an infinite volume cluster point follows
immediately. To prove relative compactness, 
we use a well-known theorem due to Prohorov.
Recall that a family $(\nu_T)$ of probability measures on 
$\Om = C([a,b],\R^d)$ is called tight if
\begin{itemize}
\item[(T1):] For all $\eta > 0$ there exists $R > 0$ such that  
 $$\nu_T(\{x \in \Om: |x_a| > R \} ) < \eta \qquad \mbox{ uniformly in } T$$.
\item[(T2):] For all $\eta>0$ and all $\varepsilon > 0$ there exists $\delta > 0$ 
such that   
  $$\nu_T(\{x \in \Om: w_{\delta}([a,b]) > \varepsilon \} ) < \eta \qquad
   \mbox{ uniformly in } T,$$
 where
  $$ w_{\delta}([a,b]) = \sup \{ |x_s-x_t| : s,t \in [a,b], |s-t| < \delta. \} $$
\end{itemize}
Prohorov's theorem states that a 
tight family of measures is relatively compact in the weak topology \cite{Bi68}. 

Usually, (T2) is rather harder to show than (T1). In our special case, however, 
(T2) follows without too much work from (T1). Without loss in generality, we 
may (and will) restrict our attention to $[a,b] = [-S,S]$.

\begin{Lemma} \label{Le3.1}
Define $(\muT)_{T > 0}$ as in (\ref{muT}), and assume (V1), (V2) and (W1). 
If $(\muT)_{T > 0}$ fulfills (T1) then it fulfills (T2) as well.
\end{Lemma}

\begin{Proof}
Fix $\eta >0$ and $\varepsilon > 0$.
Now by (T1) and the time reversibility of $\muT$ for all $T$ it is possible 
to choose $R$ such that 
$$ E_{\muT} ( |x_{-S}| > R ) < \eta/4 \quad \mbox{and} \quad 
   E_{\muT}( |x_{S}| > R ) < \eta/4 \quad \mbox{uniformly in } T.$$
Putting 
$$ B = \{ |x_{-S}| \leq R \mbox{ and } |x_{S}| \leq R \} \subset C(\R,\R^d),$$ 
and 
$$ f_{\delta}(x) = 1_{ \{w_{\delta}([-S,S]) > \varepsilon \} }(x) \quad (x \in 
C(\R,\R^d)),$$
we clearly have $\muT(B) > 1-\eta/2$ uniformly in $T$ and 
$|f_{\delta}| < 1$, and thus
\begin{equation} \label{f<1}
E_{\muT}(f_{\delta}) \leq \eta/2 + E_{\muT}(f_{\delta} 1_B) =  \eta/2 + 
E_{\muT} \left( E_{\muT} \left( f_{\delta} 1_B | \T_S \right) \right).
\end{equation}
Using Lemma \ref{fvdlr}, we find
\begin{eqnarray} \label{8Cinf}
E_{\muT} (f_{\delta} 1_B | \T_S)(\bar{x}) & = &  \frac{1}{Z_T^S(\bar{x})} 
\int e^{\cH_{\Lambda(S,T)}(x)} f_{\delta}(x) 1_B(x) \, d\mu_0^{S,\bar{x}}(x) \leq \nonumber \\
 & \leq & e^{8 C_{\infty} S} \int f_{\delta}(x) 1_B(x) \, d\mu_0^{S,\bar{x}}(x) \nonumber \\ 
 & = &  e^{8 C_{\infty} S} \int f_{\delta}(x) \, d\mu_0^{S,\bar{x}}(x) 1_B(\bar{x}).
\end{eqnarray}
The inequality above follows from (\ref{estimate}) and the definition of $Z_T^S(\bar{x})$.
Now it is easy to see that the restriction of the family 
$\{\mu_0^{S,\bar{x}}: \bar{x} \in B \}$ to $\F_{[-S,S]}$ is tight. In fact, 
this follows from the compactness of $\{x,y \in \R^d: |x| \leq R, |y| \leq R \}$.
Thus we can find 
$\delta > 0$ such that 
$$ \sup_{\bar{x} \in C(\R,\R^d)} \int f_{\delta}(x) \, d\mu_0^{S,\bar{x}}(x) 1_B(\bar{x})
    < e^{-8C_{\infty}S} \frac{\eta}{2}. $$
Using this in (\ref{8Cinf}) and plugging the resulting expression into 
(\ref{f<1}), we arrive at 
$$ \muT(f_{\delta}) \leq \eta/2 + (\eta/2) \muT(1_B) \leq \eta,$$
which is what we had to show.
\end{Proof}

\begin{Theorem} \label{Th3.2}
Assume (V1),(V2),(W1) and (W2). Then $(\mu_T)_{T > 0}$ fulfills (T1).
\end{Theorem}

\begin{Proof}
Since by (\ref{estimate}) and the stationarity of $\mu_0$ we have 
$$ e^{-2|t|C_{\infty}} \int f(x_t) \, d\muT(x) \leq \int f(x_0) \, d\muT(x) \leq 
 e^{2|t|C_{\infty}} \int f(x_t) \, d\muT(x) $$
for all $t \in \R, T > 0$ and $f \in L^{\infty}(\R^d)$, 
it will be sufficient to prove the claim for $t=0$.
We do so in several steps.\\
{ \bf Step 1:} Let $E_{\mu_0}(f|x_0=y)$ 
denote expectation with respect to the measure $\mu_0$
conditional on $x_0 = y$. Since $x \mapsto x_0$ has distribution $\psi_0^2 dx$, we have
\begin{equation} \label{3.2.1}
\muT(|x_0| > R) = \frac{1}{Z_T} \int_{|y| > R}  \psi_0^2(y) E_{\mu_0} \Big( 
e^{\cH_T} \Big| x_0=y \Big) dy.
\end{equation}
In the next few steps, we will show that there exists $K > 0$ and 
$r > 0$ such that for all $T > 0$ and all $y \in \R^d$,
\begin{equation} \label{mainestimate}
E_{\mu_0} \Big( e^{\cH_T} \Big| x_0=y \Big)
 \leq \frac{K}{\psi_0(y)} \inf_{|z| \leq r}   E_{\mu_0} \Big( 
e^{\cH_T} \Big| x_0=z \Big).
\end{equation}
Once we will have established (\ref{mainestimate}), we can plug it into 
(\ref{3.2.1}). Since moreover
$$ \frac{1}{Z_T}  \inf_{|z| \leq r}   E_{\mu_0} \Big( 
e^{\cH_T} \Big| x_0=z \Big)  
\leq \sup_{|z| \leq r} \frac{1}{\psi_0^2(z)} \muT (|x_0| \leq r) \leq \tilde{K}$$
by an expression analogous to (\ref{3.2.1}), we get
\begin{equation} \label{sdestimate}
\muT(|x_0| > R) \leq K \tilde{K} \int_{|y| > R} \psi_0(y) \, dy.
\end{equation}
The hypothesis $\psi_0 \in L^1$ from (V2) will then conclude the proof.\\
{\bf Step 2:} In order to prove (\ref{mainestimate}), we 
change the probability space we work on. Remember that $C^{(0)}$ was 
defined before equation (\ref{theta}), and consider 
\begin{equation}  \label{J}
J : C^{(0)}(\R, \R^d) \to  C([0,\infty[,\R^{2d}), \quad
        (x_t)_{t \in \R} \mapsto (x_t,x_{-t})_{t \geq 0}.
\end{equation}
$(Jx)_0 \in \R^{2d}$ is defined via the left and right hand side limits of $x_t$
as $t \to 0$, and $J$ is a bijection after making some choice for the value of 
$x \in C^{(0)}(\R,\R^d)$ at the point $0$.
We will write $\x = (x',x'')$ for the elements of $C([0,\infty[,\R^{2d})$.\\
The image of $\mu_0(.|x_0=z)$ under $J$ can be described explicitly. 
For $\z \in \R^{2d}$ denote by $\tilde{\mu}_0^{\z}$ the measure of the 
$\R^{2d}$-valued $P(\phi)_1$-process with potential
$\tilde{V}(x,y) = V(x) + V(y)$, starting in $\z$. Explicitly, if we write 
$\tilde{\F}_T$ for the $\sigma$-field over $C([0,\infty[,\R^{2d})$ generated by point 
evaluations at points within $[0,T]$, then for every $\tilde{\F}_T$-measurable, 
bounded function $f$ we have 
\begin{equation} \label{mu1}
 \int f(\x) \, d\tilde{\mu}_0^{\z}(\x) = \frac{1}{\psi_0(z')\psi_0(z'')}
 \int e^{ - \int_0^T (V(x'_s) + V(x''_s))   \, ds} f(\x) 
\psi_0(x'_T) \psi_0(x''_T) \, d\W^{\z}(\x).
\end{equation}
Here, $\W^{\z}$ denotes $2d$-dimensional Wiener measure conditional 
on $\{\x_0 = \z = (z',z'') \}$, i.e.  Brownian motion starting in $\z$.
The Markov property and time reversibility of Brownian motion together with
(\ref{FKF}) imply that
for each $z \in \R^d$,
$\tilde{\mu}_0^{(z,z)}$ is the image of $\mu_0(.|x_0=z)$ under $J$, i.e. 
$$ E_{\mu_0}(f \circ J | x_0 = z) = E_{\tilde{\mu}_0}^{(z,z)}(f).$$
Here, $E_{\tilde{\mu}_0}^{(z,z)}$ denotes expectation with respect 
to $\tilde{\mu}_0^{(z,z)}$.

Now  it is easy to check that 
\begin{eqnarray} \label{fT} 
\tilde{\cH}_T(\x) & \equiv & \cH_T \circ J^{-1}(\x) = 
- \int_0^T  \!\!\! ds \int_0^T  \!\!\! dt \,  \Big( W(x'_t,x'_s,|s-t|) + 
W(x''_t,x''_s,|s-t|) + \nonumber \\
 & &  \qquad \qquad \quad \qquad + W(x'_t,x''_s,|s+t|) + W(x''_t,x'_s,|s+t|) \Big), 
\end{eqnarray}
and therefore
\begin{equation} \label{mu1undmu0}
E_{\mu_0} \Big( e^{\cH_T} \Big| x_0=z \Big) = E_{\tilde{\mu}_0}^{(z,z)}(e^{\tilde{\cH}_T}).
\end{equation}
Thus we reduced our problem to investigating the expectation of $e^{\tilde{\cH}_T}$ 
with respect
to the strong Markov process $\tilde{\mu}_0^{\z}$ 
as a function of the starting point $\z$.\\
{\bf Step 3:} First note that in the representation established in Step 2, 
hypothesis (\ref{cond}) takes the form 
\begin{equation} \label{mainestimate2}
 \tilde{\cH}_T(\x) \leq \tilde{\cH}_T \circ \theta_{\tau} (\x) + C \tau + D 
\quad \mbox{for all } \x \in C([0,\infty[,\R^{2d}), T, \tau > 0.
\end{equation}
Here $\theta_{\tau} = J \theta^{(0)}_{\tau} J^{-1}$ is the usual time shift that maps
$(\x_t)_{t \geq 0}$ to $(\x_{t + \tau})_{t \geq 0}$. Our strategy is to use 
(\ref{mainestimate2}) together with the strong Markov property of $\tilde{\mu}_0$. 
For $r > 0$ let
$$ \tau_r(\x) = \inf \{ t \geq 0: |\x_t| \leq r \}$$
be the hitting time of the centered ball with radius $r$, and let $\F_{\tau_r}$ be
the corresponding $\sigma$-field, i.e.  
$$ \F_{\tau_r} = \{ A \in \tilde{\F}: A \cap \{\tau_r \leq t\} \in \tilde{\F}_t 
   \mbox{ for all } t \geq 0 \}. $$
Then for each $\x \in \R^{2d}$,
\begin{eqnarray} \label{StrMarkov}
 E^\x(e^{\tilde{\cH}_T}) & = &  E^\x( E^\x(e^{\tilde{\cH}_T} | \F_{\tau_r})) \leq  
 E^\x( E^\x(e^{\tilde{\cH}_T \circ \theta_{\tau_r}} e^{C \tau_r + D}| \F_{\tau_r} ) ) 
 = \nonumber \\
 & = &  E^\x(e^{C \tau_r + D}  E^\x( e^{\tilde{\cH}_T \circ \theta_{\tau_r}} | \F_{\tau_r}
  ) ) =  E^\x(e^{C \tau_r + D} E^{\x_{\tau_r}}(e^{\tilde{\cH}_T})) \leq \nonumber \\
 & \leq & \sup_{|\y| \leq r} E^\y(e^{\tilde{\cH}_T})  E^\x (e^{C \tau_r + D}).
\end{eqnarray}
All expectations above and henceforth are with respect to $\tilde{\mu}_0$. 
It remains to get a 
good estimate on the second factor on the right hand side of (\ref{StrMarkov})
and to estimate the supremum in the first factor against an infimum. This
will be done in Steps 4 and 5. \\
{\bf Step 4:} Here we show that there exists $r > 0$ and $\gamma > 0$ such that 
for all $\x \in \R^{2d}$ we have 
\begin{equation} \label{step4}
E^\x(e^{C \tau_r}) \leq 1 + \frac{C \norm[\infty]{\psi_0}}{\gamma} 
\left( \frac{1}{\psi_0(x')} + \frac{1}{\psi_0(x'')} \right) .
\end{equation}
To do so, we pick $\gamma$ with $0 < \gamma < \alpha - C$ and 
$r$ so large that $V(x) > C + \gamma$ for all
 $x \in \R^d$ with $|x| > r/ \sqrt{2}$. Obviously,
$$ \{ \x \in \R^{2d} : |\x| > r \} \subset \{ \x \in \R^{2d}: |x'| >
   r/ \sqrt{2} \} \cup  \{ \x \in \R^{2d}: |x''| >
   r/ \sqrt{2} \},$$
and with (\ref{mu1}) it follows that
\begin{eqnarray*}
\lefteqn{ \psi_0(z') \psi_0(z'') \mu_1^{\z}(\tau_r > t)  = } \\
 & = & \int e^{- \int_0^t (V(x'_s) + V(x''_s)) \, ds }
 1_{ \{|\x_s| > r \,\,\, \forall s \leq t \} } \psi_0(x'_t) 
 \psi_0(x''_t) \, d\W^\z(\x) \leq \\
 & \leq &
 \int e^{- \int_0^t V(x'_s) \, ds} e^{- \int_0^t V(x''_s) \, ds}  
  \left( 1_{ \{|(x'_s)| > r/\sqrt{2} \,\,\, \forall s \leq t \} } + 
   1_{ \{|(x''_s)| > r/\sqrt{2} \,\,\, \forall s \leq t \} } \right) \times \\
 & & \times 
 \psi_0(x'_t) \psi_0(x''_t) \, d\W^{z'}(x')
  \, d\W^{z''}(x'') = \\ 
 & = & \psi_0(z'') \int  e^{- \int_0^t V(x'_s) 
     \, ds}  1_{ \{|(x'_s)| > r/\sqrt{2} \,\,\, \forall s \leq t \} } 
   \psi_0(x'_t) \, d\W^{z'}(x') + \\
 & &  + \psi_0(z') \int  e^{- \int_0^t V(x''_s) \, ds}  
   1_{ \{|(x''_s)| > r/\sqrt{2} \,\,\, \forall s \leq t \} } 
   \psi_0(x''_t) \, d\W^{z''}(x'') \leq \\
 & \leq & (\psi_0(z') + \psi_0(z''))  
        \norm[\infty]{\psi_0} e^{-(C + \gamma) t}.
\end{eqnarray*}
The second equality above is due the eigenvalue equation
$e^{-tH_0}\psi_0 = \psi_0$ and the Feynman-Kac formula. It follows that
$$ \tilde{\mu}_0^{\z}(\tau_r > t) \leq  
   \left(  \frac{1}{\psi_0(z')} + \frac{1}{\psi_0(z'')} 
   \right) \norm[\infty]{\psi_0} e^{-(C + \gamma) t},$$
and using the equality
$$ E^\z(e^{C \tau_r}) = 1 + \int_0^{\infty} C e^{Ct} E^\z(\tau_r > t) \, dt$$
we arrive at (\ref{step4}).\\
{\bf Step 5:} Let $r > 0$ be as in Step 4. We will show that there exists
$M > 0$ such that 
 \begin{equation} \label{step5}
\sup_{|\y| \leq r} E^\y(e^{\tilde{\cH}_T}) \leq M \inf_{|\y| \leq r} E^\y(e^{\tilde{\cH}_T})
\end{equation}
uniformly in $T > 0$. Denote by $P_t(\x,\y)$ the transition density from $\x$ to 
$\y$ in time $t$ of the process $\tilde{\mu}_0$. By (\ref{mu1}) and (\ref{Transdens})
we have 
\begin{equation} \label{Transdens2}
 P_t(\x,\y) = \frac{\psi_0(y')\psi_0(y'')}{\psi_0(x') \psi_0(x'')} 
K_t(x',y') K_t(x'',y'').
\end{equation}
$\psi_0$ and $K_t$ are both uniformly bounded and bounded away from
zero on compact sets, thus for each $R > 0$ the quantity 
$$ S_t(R,r)  = \sup \left\{ \frac{P_t(\x,\z)}{P_t(\y,\z)} : \x, \y , 
\z \in \R^{2d}, |\x| \leq r, |\y | \leq r, |\z| \leq R \right \}$$
is finite. Defining $\tilde{\cH}^1_T$ like in (\ref{fT}) but with the
integrals starting at $1$ rather than at $0$,
 we see from (\ref{estimate}) that
$$  \tilde{\cH}_T(\x) - 4 C_{\infty} \leq \tilde{\cH}^1_T(\x) 
\leq \tilde{\cH}_T(\x) + 4 C_{\infty}$$
for all $\x$ and all $T$. Putting 
 $B = {\{|\x_1| < R\}}$, for each $\y$ with $|\y|< r$ we have 
\begin{equation} \label{step5.1}
E^\y(e^{\tilde{\cH}_T}) \leq  e^{4 C_{\infty}} E^\y(1_B e^{\tilde{\cH}^1_T}) + e^{C + D} 
E^\y(1_{B^c} e^{\tilde{\cH}_T \circ \theta_1}). 
\end{equation}
Defining $\bar{\cH}_T $ as in (\ref{fT}) but
with $|s+t+2|$ appearing instead of $|s+t|$ everywhere, in
the first term on the right hand side of (\ref{step5.1}) we find
\begin{eqnarray} \label{step5.2} 
 E^\y(1_B e^{\tilde{\cH}^1_T}) & = &  
 \int_{|\z| < R} P_1(\y,\z) E^\z(e^{\bar{\cH}_{T-1}}) \, d\z \leq \nonumber \\
 & \leq & S_1(R,r) \int_{|\z| \leq R} P_1(\x,\z) E^\z(e^{\bar{\cH}_{T-1}}) 
 \, d\z = \nonumber \\
 & = & S_1(R,r) E^\x(1_B e^{\tilde{\cH}^1_T}) \leq S_1(R,r)  e^{4 C_{\infty}} 
       E^\x(e^{\tilde{\cH}_T})
\end{eqnarray}
for each $\x$ with $|\x| \leq r$. Turning to the second term on the right hand side 
of (\ref{step5.1}), equations (\ref{StrMarkov}) and (\ref{step4}) give 
\begin{eqnarray} \label{step5.3}
\lefteqn{ E^\y(1_{B^c} e^{\tilde{\cH}_T \circ \theta_1}) = 
   \int_{|\z| > R} P_1(\y,\z) E^\z(e^{\tilde{\cH}_T}) \, d\z  \leq \nonumber } \\
 & \leq & \sup_{|\x| \leq r} E^{\x}(e^{\tilde{\cH}_T}) \int_{|\z| > R} P_1(\y,\z) 
    E^\z(e^{C \tau_r + D})
 \, d\z \leq  \\
 & \leq & \sup_{|\x| \leq r} E^{\x}(e^{\tilde{\cH}_T}) e^D  \int_{|\z| > R} P_1(\y,\z)
    \left( 1 + \frac{C \norm[\infty]{\psi_0}}{\gamma} 
    \left( \frac{1}{\psi_0(z')} + \frac{1}{\psi_0(z'')} \right) \right) 
    d\z. \nonumber 
\end{eqnarray}
By (\ref{Transdens2}) and the eigenvalue equation, we have
\begin{eqnarray} \label{Sternchen}
\lefteqn{ \int P_1(\y,\z)  
 \left( \frac{1}{\psi_0(z')} + \frac{1}{\psi_0(z'')} \right) \, d\z = } \nonumber \\
& = & \frac{1}{\psi_0(y')} \int K_1(y'',z) \, dz +  
       \frac{1}{\psi_0(y'')} \int  K_1(y',z) \, dz.
\end{eqnarray}
By (V1), the above integrals are bounded in $y'$ and $y''$, respectively \cite{Si82},
and thus the right hand side of (\ref{Sternchen}) is 
uniformly bounded on 
$\{ \y: |\y| < r \}$. This implies that there exists $\bar{R} > 0$ and $\delta < 1$ 
such that 
$$  \int_{|\z| > \bar{R}} P_1(\y,\z)
    \left( 1 + \frac{C \norm[\infty]{\psi_0}}{\gamma} 
    \left( \frac{1}{\psi_0(z')} + \frac{1}{\psi_0(z'')} \right) \right) 
    d\z \leq e^{-(C + 2D)} \delta $$ 
uniformly on $\{ \y: |\y| < r \}$. 
Plugging this result together with (\ref{step5.2})
into (\ref{step5.1}), we arrive at
\begin{equation} \label{step5.4}
 E^\y(e^{\tilde{\cH}_T}) \leq  S_1(\bar{R},r)  e^{8 C_{\infty}} E^\x(e^{\tilde{\cH}_T}) +
 \delta  \sup_{|\z| \leq r} E^{\z}(e^{\tilde{\cH}_T}),
\end{equation}
which is valid for all $\x,\y$ with $|\x|,|\y| \leq r$. By taking the supremum over $\y$ 
and the infimum over $\x$ in (\ref{step5.4}) and rearranging, we find
$$ \sup_{|\y| \leq r} E^{\y}(e^{\tilde{\cH}_T}) \leq  \frac{S_1(\bar{R},r)  
e^{8 C_{\infty}}}{1 - \delta}  \inf_{|\y| \leq r} E^{\y}(e^{\tilde{\cH}_T}),$$
which concludes Step 5 and the proof.
\end{Proof} 

The two previous statements show relative compactness of the restrictions
$\{\mu_T|_{\F_{[-S,S]}}: T > 0 \}$ for any $S > 0$. From here, it is only a small
step to relative compactness in the topology of local weak convergence.

\begin{Theorem} \label{Th3.3}
Assume (V1),(V2),(W1) and (W2). Then $\{\muT : T \geq 0 \}$ is relatively compact
in the topology of local weak convergence. Consequently, the family has an
infinite volume cluster point.
\end{Theorem}

\begin{Proof}
Take $S > 0$ and fix any sequence $(T_n) \subset \R^+$. 
By Lemma \ref{Le3.1}, Theorem \ref{Th3.2} and the tightness argument, for each fixed $S>0$
there exists a subsequence $(t_n)$ of $(T_n)$ such that
$(\mu_{t_n}|_{\F_{[-S,S]}})$ converges weakly to some probability measure $\mu_{\infty}$ on
$C([-S,S], \R^d)$ . 
In case $L = \limsup_{n \to \infty} T_n < \infty$ we are done
by choosing $S > L$. In case $L = \infty$, we observe that convergence of 
$(\mu_{t_n}|_{\F_{[-R,R]}})_{n \in \N}$ implies convergence of 
$(\mu_{t_n}|_{\F_{[-S,S]}})_{n \in \N}$ if $R > S$, and thus
a diagonal sequence argument does the job. This second case also provides us
with an infinite volume cluster point.
\end{Proof}

Let us denote by $\mu$ any cluster point of the family $(\muT)_{T > 0}$ obtained by
Theorem \ref{Th3.3}.  
Due to the good control on the stationary density we obtain in Theorem \ref{Th3.2}, 
we have the following estimate on the growth of paths under $\mu$.

\begin{Lemma} \label{Le4.1}
Let $f:\R^+ \to \R^+$ be monotone increasing with $f(x) \to \infty$ as $x \to \infty$, 
and suppose that
\begin{equation} \label{pathcond}
 \sum_{n=1}^{\infty} \int_{|y| > f(n)} \psi_0(y) \, dy < \infty
\end{equation}
Then for $\mu$-almost every path $x \in C(\R, \R^d)$, we have
$$ \limsup_{|t| \to \infty} \frac{|x_t|}{f(|t|)} \leq 1.$$
\end{Lemma}

\begin{Proof}
By path continuity and time reversibility, it is obviously enough to prove that for each 
$k \in \N$,
$$ \mu \left( \limsup_{n \to \infty} \frac{|x_{(n/k)}|}{f(n/k)} > 1 \right) = 0.$$ 
Since the above event is equal to $\{x: |x_{n/k}| > f(n/k) \mbox{ infinitely often} \}$,
the first Borel-Cantelli lemma will yield the result once we have checked that 
\begin{equation} \label{bccond}
\sum_{n=1}^{\infty} \mu( |x_{n/k}| > f(n/k) ) < \infty.
\end{equation}
By the stationarity of $\mu$ and equation (\ref{sdestimate}),
there exists a constant $M$ such that  
$$ \mu( |x_{n/k}| > f(n/k) ) = \mu( |x_0| > f(n/k) ) \leq M
   \int_{|y| > f(n/k)} \psi_0(y) \, dy $$
for $n$ large enough. Since
$$ \sum_{n=1}^{\infty} \int_{|y| > f(n/k)} \psi_0(y) \, dy \leq k \sum_{n=1}^{\infty}
  \int_{|y| > f(n)} \psi_0(y) \, dy,$$
(\ref{pathcond}) implies (\ref{bccond}).
\end{Proof}

In many cases, estimates on the decay of $\psi_0$ can be obtained via $V$. In \cite{Ca78}
it is shown that for $s \geq 0$ the estimate
$\liminf_{|x| \to \infty} V(x)/|x|^{2s} > 0$ implies the existence of constants $A>0,\beta>0$ 
such that 
$$ \psi_0(y) \leq A \exp(-\beta|y|^{s+1})$$
for all $y \in \R^d$. In this case,  Lemma \ref{Le4.1} implies
$$ \limsup_{|t| \to \infty} \frac{|x_t|}{(\gamma \ln{|t|})^{s+1}} = 0$$
for each $\gamma > 1/\beta$
and $\mu$-almost all $x \in C(\R,\R^d)$. This result has been obtained (for $s>1$) in
\cite{LM00} via the cluster expansion.

We conclude this paper by showing the 
infinite volume analogue of Lemma \ref{fvdlr}. We refer to 
Section \ref{S2} for notation and additionally introduce
\begin{eqnarray}
\Lambda(S) & = & (\R \times [-S,S]) \cup ([-S,S] \times \R), \label{LambdaS}\\
d\mu^{S,\bar{x}} & = & 
 \frac{1}{Z^S(\bar{x})} \exp(\cH_{\Lambda(S)}(x)) \, d\mu_0^{S,\bar{x}}(x). \label{infvol}
\end{eqnarray}
Note that the normalizing constant $Z^S(\bar{x})$ is finite for each 
$\bar{x} \in C(\R,\R^d)$ due to (\ref{estimate}).

\begin{Proposition} \label{infvoldlr}
For each $S > 0$ and each infinite volume cluster point $\mu$ of $(\mu_T)$, 
$\bar{x} \mapsto \mu_T^{S,\bar{x}}$ is a version of the regular 
conditional probability $\mu(.| \T_S)$. In other words, $\mu$ is a Gibbs
measure for the reference measure $\mu_0$ and the potential $W$.
\end{Proposition}

\begin{Proof}
By Lemma \ref{fvdlr}, we have for $f,g \in L^{\infty}(C(\R,\R^d))$ with 
$\T_S$-measurable $g$ that
\begin{equation} \label{fivers}
 \int g(\bar{x}) E_{\muT^{S,\bar{x}}}(f) \, d\muT(\bar{x}) = E_{\muT}(fg).
\end{equation}
We have to show that (\ref{fivers}) remains true when we replace $\muT$ by $\mu$
and $\muT^{S,\bar{x}}$ by $\mu^{S,\bar{x}}$.
By a monotone class argument, we may in assume that $f$ and $g$ are 
$\F_{[-R,R]}$-measurable for some $R>S$. Taking a sequence $(t_n)$ such that $\mu_{t_n}$
converges to $\mu$, we immediately see that the right hand side of (\ref{fivers})
converges to $E_{\mu}(fg)$. As for the left hand side, (\ref{estimate}) guarantees
that $\cH_{\Lambda(S,T)}(q)$ converges to $\cH_{\Lambda(S)}(q)$ uniformly in 
$q \in C(\R,\R^d)$ as $T \to \infty$, and thus the left hand side converges to 
$\int g(\bar{x}) E_{\mu^{S,\bar{x}}}(f) \, d\mu(\bar{x})$.
\end{Proof}

{\bf Acknowledgment:} I wish to thank Prof. H. Spohn for constant encouragement and
uncountably many useful discussions.

\end{document}